\documentclass{article}
\usepackage{amsfonts}
\usepackage{amsmath}

\newtheorem{theorem}{Theorem}

\newtheorem{corollary}[theorem]{Corollary}

\newtheorem{example}[theorem]{Example}

\newtheorem{lemma}[theorem]{Lemma}

\newtheorem{proposition}[theorem]{Proposition}
\newtheorem{remark}[theorem]{Remark}

\newenvironment{proof}[1][Proof]{\noindent\textbf{#1.} }{\ \rule{0.5em}{0.5em}}

\begin{document}

\title{A generalization of Osgood's test and a comparison criterion for
integral equations\\
with noise}
\author{M. J. Ceballos-Lira \\
Divisi\'{o}n Acad\'{e}mica de Ciencias B\'{a}sicas,\\
Universidad Ju\'{a}rez Aut\'{o}noma de Tabasco,\\
Km. 1 Carretera Cunduac\'{a}n-Jalpa de M\'{e}ndez, \\
Cunduac\'{a}n, Tab. 86690, Mexico \and J. E. Mac\'{\i}as-D\'{\i}az and J.
Villa \\
Departamento de Matem\'{a}ticas y F\'{\i}sica, \\
Universidad Aut\'{o}noma de Aguascalientes,\\
Avenida Universidad 940, Ciudad Universitaria, \\
Aguascalientes, Ags. 20131, Mexico\\
\texttt{jvilla@correo.uaa.mx}}
\date{}
\maketitle

\begin{abstract}
In this work, we prove a generalization of Osgood's test for the explosion
of the solutions of initial-value problems. We also establish a comparison
criterion for the solution of integral equations with noise, and provide
estimations of the time of explosion of problems arising in the
investigation of crack failures where the noise is the absolute value of the
Brownian motion.
\end{abstract}

{\small {\noindent }\textbf{Key words:} Osgood's test, comparison criterion,
time of explosion, integral equations with noise, crack failure}


\section{Introduction\label{S:Intro}}

Let $x_{0}$ be a positive, real number, let $b$ be a positive, real-valued
function defined on $[0,\infty )$, and suppose that $y$ is an
extended-real-valued function with the same domain as $b$. The present work
is motivated by a criterion for the explosion of the solutions of ordinary
differential equations of the form 
\begin{equation}
\left\{ 
\begin{array}{ll}
\displaystyle{\frac{dy(t)}{dt}=b(y(t))}, & t>0, \\ 
y(0)=x_{0}. & 
\end{array}
\right.  \label{Eq:ODE}
\end{equation}
More precisely, the time of explosion of the solution of this initial-value
problem is the nonnegative, extended-real number $t_{e}=\sup \{t\geq
0:y(t)<\infty \}$. The above-mentioned criterion is called \emph{Osgood's
test} after its author \cite{Osgood}, and it states that $t_{e}$ is finite
if and only if $\int_{x_{0}}^{\infty }ds/b(s)<\infty $. In such case, $%
t_{e}=\int_{x_{0}}^{\infty }ds/b(s)$.

A natural question readily arises about the possibility to extend Osgood's
test to more general, initial-value problems, say, to problems in which the
drift function $b$ in the ordinary differential equation of (\ref{Eq:ODE})
is multiplied by a suitable, nonnegative function of $t$. Another direction
of investigation would be to investigate conditions under which the
solutions of the integral form of such equation with a noise function added,
explode in a finite time (see for example \cite{L-V}). Evidently the
consideration of these two problems as a single one is an interesting topic
of study \emph{per se}. In fact, the purpose of this paper is to provide a
generalization of Osgood's test to integral equations with noise, which
generalize the problem presented in (\ref{Eq:ODE}). Important, as it is in
the recent literature \cite{Du, Faramarz, Kafini, Malolepszy}, the problem
of establishing analytical conditions under which the time of explosion of
the problem under investigation, is tackled here. In the way, we establish a
comparison criterion for the solutions of integral equations with noise, and
show some applications to the spread of cracks in rigid surfaces.

Our manuscript is divided in the following way: Section \ref{S:Osgood}
introduces the integral equation with noise that motivates this manuscript,
along with a convenient simplification for its study; a generalization of
Osgood's test is presented in this stage for the associated initial-value
problem for both scenarios: noiseless and noisy systems. Section \ref%
{S:Compar} establishes a comparison criterion for the solutions of two
noiseless systems with comparable initial conditions. A necessary condition
for the explosion of the solutions of the problem under investigation is
provided in this section, together with an illustrative counterexample and a
partial converse. In Section \ref{S:Approx}, we give upper and lower bounds
for the value of the time of explosion of our integral equation. Finally,
Section \ref{S:Applic} provides estimates of probabilities associated to the
time of explosion of a system in which the noise is the absolute value of
the Brownian motion.

\section{Osgood's test\label{S:Osgood}}

Let $\overline {\mathbb{R}}$ denote the set of extended-real numbers.
Throughout, $a, b : [0 , \infty) \rightarrow \mathbb{R}$ will represent
positive, continuous functions, while the function $g : [0 , \infty)
\rightarrow \mathbb{R}$ will be continuous and nonnegative. For physical
reasons, the function $g$ is called a \emph{noise}. In this work, $x _0$
will denote a positive, real number, and $X : [0 , \infty) \rightarrow 
\overline {\mathbb{R}}$ will be a nonnegative function whose dependency on $%
t \geq 0$ is represented by $X _t$. We are interested in establishing
conditions under which the solutions of the integral equation 
\begin{equation}  \label{Eq:Model}
X _t = x _0 + \int _0 ^t a (s) b (X _s) d s + g (t), \quad t \geq 0,
\end{equation}
explode in finite time. More precisely, we define the \emph{time of explosion%
} of $X$ as the nonnegative, extended-real number $T ^X _e = \sup \{t \geq 0
: X _t < \infty \}$. In this manuscript, we investigate conditions under
which the time of explosion of $X$ is a real number.

Letting $Y _t = X _t - g (t)$, one sees immediately that the problem under
consideration is equivalent to finding the time of explosion of the solution 
$Y$ of the equation 
\begin{equation}  \label{Eq:EIM}
Y _t = x _0 + \int _0 ^t a (s) b (Y _s + g (s)) d s, \quad t \geq 0.
\end{equation}
As a matter of fact, $T _e ^X = T _e ^Y$. From this point on, this common,
extended-real number will be denoted simply by $T _e$ for the sake of
briefness.

\begin{remark}
\emph{It is worth noticing that (\ref{Eq:EIM}) can be presented in
differential form as the equivalent, initial-value problem 
\begin{equation}
\left\{ 
\begin{array}{ll}
\displaystyle{\frac{dY_{t}}{dt}=a(t)b(Y_{t}+g(t))}, & t>0, \\ 
Y_{0}=x_{0}, & 
\end{array}%
\right.  \label{Eq:EDM}
\end{equation}%
a problem for which the existence of solutions is guaranteed, for instance,
when $b$ is locally Lipschitzian and $a$ is regulated (see (10.4.6) in \cite%
{Dieudonne}).}
\end{remark}

Let $r$ be a real number such that $0<r\leq x_{0}$. We define the functions $%
A:[0,\infty )\rightarrow \mathbb{R}$ and $B_{r}:[x_{0}-r,\infty )\rightarrow 
\mathbb{R}$ by 
\begin{equation*}
A(t)=\int_{0}^{t}a(s)ds\quad \text{and}\quad B_{r}(x)=\int_{x_{0}-r}^{x}%
\frac{ds}{b(s)}.
\end{equation*}%
For the sake of convenience, we let $B$ be the function $B_{0}$. Evidently,
both of these functions are nonnegative, increasing and continuous, and so
are their inverses. On the other hand, if $r\geq -x_{0}$, we let $\widetilde{%
B}_{r}:[x_{0}+r,\infty )\rightarrow \mathbb{R}$ be given by $\widetilde{B}%
_{r}=B_{-r}$. For every fixed $x\geq x_{0}$, we define $\widetilde{B}%
^{x}:[-x_{0},x-x_{0}]\rightarrow \mathbb{R}$ by $\widetilde{B}^{x}(r)=%
\widetilde{B}_{r}(x)$; we prefer this second notation in either case.
Additionally, we define $\beta :[-x_{0},\infty )\rightarrow \mathbb{R}$ by $%
\beta (t)=\widetilde{B}^{\infty }(t)$. All of these functions and their
inverses are nonnegative, continuous and decreasing in their domains.

\begin{lemma}[Generalized Osgood's test]
The initial-value problem 
\begin{equation}  \label{Eq:EGO}
\left\{ 
\begin{array}{ll}
\displaystyle {\frac {d y (t)} {d t} = a (t) b (y (t))}, & t > 0, \\ 
y (0) = x _0, & 
\end{array}
\right.
\end{equation}
has a unique solution given by $y (t) = B ^{- 1} (A (t))$, for $t < A ^{- 1}
(B (\infty))$. The solution explodes in finite time if and only if $B
(\infty) < A (\infty)$, in which case, $T _e ^y = A ^{- 1} (B (\infty))$.
\end{lemma}

\begin{proof}
The function $y(t)=B^{-1}(A(t))$ is evidently a solution of (\ref{Eq:EGO}).
Additionally, expressing the differential equation in (\ref{Eq:EGO}) as $%
y^{\prime }(s)/b(y(s))=a(s)$, integrating both sides over $[0,t]$ and
performing a suitable substitution, we obtain that $B(y(t))=A(t)$, whence
the uniqueness follows. Moreover, $y(t)$ is real if and only if $%
t<A^{-1}(B(\infty ))$.

Now, if the solution of (\ref{Eq:EGO}) explodes at the time $t_{e}<\infty $,
then $B(\infty )=A(t_{e})<A(\infty )$. Conversely, the number $%
A^{-1}(B(\infty ))$ is real, so that 
\begin{equation*}
B(y(A^{-1}(B(\infty ))))\allowbreak =\allowbreak A(A^{-1}(B(\infty
)))\allowbreak =\allowbreak B(\infty ).
\end{equation*}
This implies that $T_{e}^{y}\leq A^{-1}(B(\infty ))$, and the opposite
inequality follows from the fact that the solution of (\ref{Eq:EGO}) exists
for $t<A^{-1}(B(\infty ))$.\hfill
\end{proof}

As a consequence, the solution of (\ref{Eq:EGO}) is nonnegative, continuous
and increasing on $[0,T_{e}^{y})$, and so is its inverse on $[x_{0},\infty )$%
. Likewise, the function $\overline{B}:[x_{0},\infty )\rightarrow \mathbb{R}$%
, given by the formula 
\begin{equation*}
\overline{B}(y)=\int_{x_{0}}^{y}\frac{ds}{b(s+g(Y^{-1}(s)))},
\end{equation*}%
is nonnegative, continuous and increasing.

\begin{corollary}
\label{Coro:1} The solution of (\ref{Eq:EDM}) can be expressed as $Y_{t}=%
\overline{B}\,^{-1}(A(t))$, for every $t<A^{-1}(\overline{B}(\infty ))$.
\end{corollary}

\begin{proof}
The proof runs as in Osgood's test. \hfill
\end{proof}

\section{A comparison theorem\label{S:Compar}}

\begin{theorem}[Comparison criterion]
\label{Thm:Comparison} Let $0<x_{0}\leq x_{1}$, let $b$ be
non-de\-crea\-sing, and assume that the functions $u,v:[0,\infty
)\rightarrow \overline{\mathbb{R}}$ satisfy 
\begin{equation*}
v(t)\geq x_{1}+\int_{0}^{t}a(s)b(v(s))ds\quad \text{and}\quad
u(t)=x_{0}+\int_{0}^{t}a(s)b(u(s))ds,\quad t\geq 0.
\end{equation*}%
Then, $v(t)\geq u(t)$ for every $t\geq 0$, and $T_{e}^{v}\leq
A^{-1}(B(\infty ))$.
\end{theorem}

\begin{proof}
It is sufficient to show that $v\geq u$ because, in such case, $%
T_{e}^{v}\leq T_{e}^{u}=A^{-1}(B(\infty ))$. Assume first that $x_{0}<x_{1}$%
, and let $N=\{t\geq 0:u(s)\leq v(s),s\in \lbrack 0,t]\}$. The set $N$ is
nonempty, so $\widetilde{T}=\sup N$ exists in $\overline{\mathbb{R}}$. If $%
\widetilde{T}$ were a real number, then 
\begin{equation*}
L=\lim_{\epsilon \rightarrow 0^{+}}(v(\widetilde{T}+\epsilon )-u(\widetilde{T%
}+\epsilon ))\geq x_{1}-x_{0}+\lim_{\epsilon \rightarrow 0^{+}}\int_{%
\widetilde{T}}^{\widetilde{T}+\epsilon }a(s)\left[ b(v(s))-b(u(s))\right] ds
\end{equation*}%
by the fact that $v(s)-u(s)\geq 0$, for every $s\in \lbrack 0,\widetilde{T}]$%
. It follows that $L\geq x_{1}-x_{0}$. By definition, there exists $\delta
>0 $ such that $v(\widetilde{T}+s)-u(\widetilde{T}+s)>0$ for every $s\in
\lbrack 0,\delta )$, whence it follows that $\widetilde{T}+\frac{\delta }{2}%
\in N$, a contradiction. Consequently, $u(t)\leq v(t)$ for every $t\geq 0$.
Now, in case that $x_{0}=x_{1}$, the solution of the equation 
\begin{equation*}
u_{r}(t)=x_{0}-r+\int_{0}^{t}a(s)b(u_{r}(s))ds,\quad 0<r<x_{0},
\end{equation*}%
satisfies $v(t)\geq u_{r}(t)$, for every $t\geq 0$. Using Osgood's test and
the continuity of $B_{r}^{-1}$, we obtain 
\begin{equation*}
v(t)\geq \lim_{r\rightarrow 0^{+}}u_{r}(t)=\lim_{r\rightarrow
0^{+}}B_{r}^{-1}(A(t))=B^{-1}(A(t))=u(t).
\end{equation*}
\hfill
\end{proof}

\begin{theorem}
\label{Thm:2} Suppose that $b$ is non-decreasing, and $B(\infty )<A(\infty )$%
. Then the solution of (\eqref {Eq:EDM}) explodes in finite time. The time
of explosion of $Y$ is $t_{e}=A^{-1}\left( \overline{B}(\infty )\right) $.
\end{theorem}

\begin{proof}
The fact that $b$ is non-decreasing yields 
\begin{equation*}
Y_{t}=x_{0}+\int_{0}^{t}a(s)b(Y_{s}+g(s))ds\geq
x_{0}+\int_{0}^{t}a(s)b(Y_{s})ds.
\end{equation*}%
Theorem \ref{Thm:Comparison} gives that $T_{e}\leq A^{-1}(B(\infty ))$ when
we compare $Y$ with the solution of $\widetilde{Y}_{t}=x_{0}+%
\int_{0}^{t}a(s)b(\widetilde{Y}_{s})ds$. On the other hand, $A(t_{e})=%
\overline{B}(\infty )\leq B(\infty )<A(\infty )$, which implies that $t_{e}$
is real. The expression of $t_{e}$ and Corollary \ref{Coro:1} yield $%
A(t_{e})=\overline{B}(Y_{T_{e}})=A(T_{e})$. We conclude that $T_{e}=t_{e}$.\hfill
\end{proof}

Intuitively, it is not generally true that the explosion of the solutions of
(\ref{Eq:EDM}) in finite time is a sufficient condition for the inequality $%
B(\infty )<A(\infty )$ to be satisfied. This assertion follows after
noticing that $A$ and $B$ do not depend of the noise function; however, the
time of explosion does. We will establish our claim precisely through the
following counter-example.

\begin{example}
\emph{Let $x_{0}=1$, and let $a$, $b$ and $g$ be given by the expressions $%
a(t)=e^{-t}$, $b(t)=\frac{1}{4}t^{3}$, and $g(t)=e^{t}$, for every $t>0$.
Expanding the expression $(Y_{s}+e^{s})^{3}$ in (\ref{Eq:EIM}), we obtain $%
Y_{t}\geq 1+\frac{1}{4}\int_{0}^{t}Y_{s}^{2}ds$. Then $Y_{t}\geq (1-\frac{1}{%
4}t)^{-1}$, which implies that $Y$ explodes in finite time. However, $%
B(\infty )=2>1=A(\infty )$.}
\end{example}

The following result is a partial converse of Theorem \ref{Thm:2}. We let $%
\widehat {g} (t) = \sup \{ g (s) : s \in [0 , t] \}$, for every $t \geq 0$.

\begin{proposition}
\label{Prop:1} Suppose that $b$ is non-decreasing, and that 
\begin{equation*}
\widehat{g}(t)<b(x_{0})\int_{t}^{\infty }a(s)ds.
\end{equation*}%
If the solution $Y$ of (\ref{Eq:EIM}) explodes in finite time, then $%
B(\infty )<A(\infty )$.
\end{proposition}

\begin{proof}
By Corollary \ref{Coro:1}, $\overline{B}(\infty )=\overline{B}%
(Y_{T_{e}})=A(T_{e})$. Since $b$ is non-decreasing and $g(Y_{s}^{-1})\leq 
\widehat{g}(T_{e})$ for every $s\in \lbrack x_{0},\infty )$, we obtain 
\begin{equation*}
\int_{x_{0}}^{\infty }\frac{ds}{b(s+\widehat{g}(T_{e}))}\leq A(T_{e}).
\end{equation*}%
Separating the integral in the definition of $B(\infty )$ as the sum of the
integrals over the intervals $[x_{0},x_{0}+\widehat{g}(T_{e})]$ and $[x_{0}+%
\widehat{g}(T_{e}),\infty )$, using the facts that $b$ is positive and
non-decreasing, and employing the last inequality and the hypothesis, we
obtain 
\begin{eqnarray*}
B(\infty ) &\leq &\int_{x_{0}}^{x_{0}+\widehat{g}(T_{e})}\frac{ds}{b(x_{0})}%
+\int_{x_{0}}^{\infty }\frac{ds}{b(s+\widehat{g}(T_{e}))} \\
&\leq &\frac{\widehat{g}(T_{e})}{b(x_{0})}+A(T_{e})<A(\infty ).
\end{eqnarray*}
\hfill
\end{proof}

\section{Approximation of the explosion time\label{S:Approx}}

It is important to notice that the time of explosion of $Y$, as given by the
Theorem \ref{Thm:2}, presents the disadvantage of depending on the solution $%
Y$ itself. In this section, we will derive some approximations to $T _e$
which do not present this shortcoming. For the remainder of this manuscript
and for the sake of convenience, we let $T = A ^{- 1} (B (\infty))$.
Throughout this section, $b$ will be a non-decreasing function.

The Comparison criterion and Osgood's test yield that the time of explosion
of the solution $Y$ of (\ref{Eq:EIM}) satisfies $T_{e}\leq T$. On the other
hand, 
\begin{equation*}
Y_{t}\leq x_{0}+\int_{0}^{t}a(s)b(Y_{s}+\widehat{g}(T))ds,
\end{equation*}%
and the Comparison criterion leads us to conclude that 
\begin{equation}
A^{-1}\left( \beta (\widehat{g}(T))\right) \leq T_{e}\leq T.  \label{Eq:EPA}
\end{equation}

In general, the function $b : [0 , \infty) \rightarrow \mathbb{R}$ is \emph{%
sub-multiplicative} if there exists a positive constant $c$ such that $b (x
y) \leq c b (x) b (y)$, for every $x , y \geq 0$. Evidently, exponential and
power functions are sub-multiplicative.

Suppose that $b$ is a sub-multiplicative function, and let $c$ be the
positive number provided by the definition of sub-multiplicativity. In the
following, it will be convenient to define the function $\widetilde{A}%
:[0,\infty )\rightarrow \mathbb{R}$ by 
\begin{equation*}
\widetilde{A}(t)=c\int_{0}^{t}a(s)b\left( \frac{1}{x_{0}}g(s)+1\right) ds.
\end{equation*}%
This function is nonnegative, continuous and increasing and, thus, it is
invertible, and has a continuous and increasing inverse.

\begin{proposition}
\label{Prop:2} Let $b$ be a sub-multiplicative function. Then $T _e \geq 
\widetilde {A} ^{- 1} (B (\infty))$.
\end{proposition}

\begin{proof}
Since $Y_{s}\geq x_{0}$ for every $s\geq 0$, we obtain 
\begin{eqnarray*}
g(s)+Y_{s} &=&Y_{s}g(s)\left( \frac{1}{Y_{s}}+\frac{1}{g(s)}\right) \\
&\leq &Y_{s}g(s)\left( \frac{1}{x_{0}}+\frac{1}{g(s)}\right) =Y_{s}\left( 
\frac{1}{x_{0}}g(s)+1\right) .
\end{eqnarray*}%
Monotonicity and sub-multiplicativity of $b$, along with \eqref {Eq:EIM},
yield 
\begin{equation*}
Y_{t}\leq x_{0}+c\int_{0}^{t}a(s)b\left( \frac{1}{x_{0}}g(s)+1\right)
b(Y_{s})ds.
\end{equation*}%
The conclusion of the proposition follows now from the Comparison criterion
and Osgood's test. \hfill
\end{proof}

\section{An application\label{S:Applic}}

Throughout this section, we consider the stochastic differential equation (%
\ref{Eq:Model}) with noise function $|W_{t}|$, where $W$ is the Brownian
motion. The noise function is taken in absolute value in view of physical
considerations on the dynamics of cracks growth under fatigue loading. In
fact, it has been established experimentally that cracks in the subcritical
stage grow with a velocity that increases with the crack length \cite{Paris}%
. The governing equation is called \emph{Paris' equation}, and it is a power
law (which is a sub-multiplicative function) in which the exponent is
determined empirically. As a matter of fact, it has been established that
Paris' law is valid for a wide range of materials \cite{Allen1, Allen2,
Mach, Wei}.

For the remainder of this work, we let $\Phi (x)$ represent the probability
that a random variable with standard normal distribution assumes values in $%
[0 , x]$, for every $x \geq 0$.

\begin{proposition}
Let $0 \leq t < T$. Then 
\begin{equation}  \label{Eq:PRA1}
P (T _e \leq t) \leq 1 - \Phi \left( \frac {\beta ^{- 1} (A (t))} {\sqrt {T}}
\right).
\end{equation}
\end{proposition}

\begin{proof}
We use here the first inequality of (\ref{Eq:EPA}). Notice that 
\begin{equation*}
P\left( T_{e}\leq t\right) \leq P\left( A^{-1}(\beta (|\widehat{W}%
_{T}|))\leq t\right) =P\left( |\widehat{W}_{T}|\geq \beta ^{-1}(A(t))\right)
.
\end{equation*}
Equation (8.4) in \cite{Karatzas} completes the proof. \hfill
\end{proof}

For every nonnegative, real number $r$, we let $T _r = \inf \{ t > 0 : | W
_t | = r \}$. Evidently, $| W _s | \leq r$, for every $s \in [0 , T _r]$.

\begin{proposition}
\label{Prop:4} Let $0 \leq t \leq T$. For every $r \geq 0$, 
\begin{equation}  \label{Eq:PRA2}
P \left( T _e \leq t | T _r < T \right) \leq \frac {1 - \Phi \left(r /\sqrt {%
A ^{- 1} (B (\widetilde {B} ^{- 1} _r (A (t))))}\right)} {1 - \Phi \left(r / 
\sqrt {T}\right)}.
\end{equation}
\end{proposition}

\begin{proof}
Notice that $|\widehat{W}_{T}|\geq r$ whenever $T_{r}<T$. Moreover, Osgood's
test and the Comparison criterion imply that $Y_{t}\geq B^{-1}(A(t))$, for
every $t\geq 0$. Using (\ref{Eq:EPA}), we obtain $A(T_{e})\geq \widetilde{B}%
_{|\widehat{W}_{T}|}(\infty )\geq \widetilde{B}_{r}(\infty )\geq \widetilde{B%
}_{r}(B^{-1}(A(T_{r})))$. Therefore, 
\begin{eqnarray*}
P\left( T_{e}\leq t|T_{r}<T\right) &\leq &\frac{P\left( \widetilde{B}%
_{r}(B^{-1}(A(T_{r})))\leq A(t)\right) }{P(T_{r}<T)} \\
&=&\frac{P\left( T_{r}<A^{-1}(B(\widetilde{B}_{r}^{-1}(A(t))))\right) }{%
P(T_{r}<T)}.
\end{eqnarray*}%
The conclusion follows now from \cite{Karatzas} as in Proposition \ref%
{Prop:4}. \hfill
\end{proof}

On physical grounds, the function $Y$ may represent the temporal behavior of
the transversal length of a crack failure on some material \cite{Paris}. In
this context, the parameter $x _0$ represents the initial, transversal
length of the crack, and $L$ is the transversal length of the material. For
practical purposes, one may think of the wing of an airplane which has a
fixed transversal length, on which there is a crack with known initial
length. In such case, one investigates the dynamics of the length of the
crack with respect to time, in order to conduct preventive maintenance on
the wing and avoid possible accidents \cite{Augustin}.

\begin{proposition}
If $L>x_{0}$, then 
\begin{equation*}
P\left( Y_{L}^{-1}\leq t\right) \leq 1-\Phi \left( \frac{\widetilde{B}%
_{A(t)}^{-1}(L)}{\sqrt{T}}\right) .
\end{equation*}
\end{proposition}

\begin{proof}
Let $\widetilde{Y}$ the solution of $\widetilde{Y}_{t}=x_{0}+%
\int_{0}^{t}a(s)b(\widetilde{Y}_{s}+|\widehat{W}_{T}|)ds$, for every $0\leq
t<T$. By Osgood's test and the Comparison criterion, $\widetilde{B}_{|%
\widehat{W}_{T}|}^{-1}(A(t))=\widetilde{Y}_{t}\geq Y_{t}$. Once again, the
conclusion is reached using \cite{Karatzas} in the right-most end of the
chain of identities and inequalities 
\begin{eqnarray*}
P\left( Y_{L}^{-1}\leq t\right) &\leq &P\left( \widetilde{Y}_{L}^{-1}\leq
t\right) \\
&=&P\left( \widetilde{B}_{|\widehat{W}_{T}|}(L)\leq A(t)\right) \\
&=&1-P\left( |\widehat{W}_{T}|\leq \widetilde{B}_{A(t)}^{-1}(L)\right) .
\end{eqnarray*}
\hfill
\end{proof}

\begin{example}
\emph{Let $x_{0}$, $a_{0}$ and $\alpha $ be positive numbers, and let $%
a(t)=a_{0}$ and $b(t)=t^{1+\alpha }$, for every $t\geq 0$. Observe that $%
A(t)=a_{0}t$ and, for every $r\in \lbrack 0,x_{0}]$ and every $x\geq x_{0}$, 
\begin{equation*}
B_{r}(x)=\frac{1}{\alpha }\left[ \frac{1}{(x_{0}-r)^{\alpha }}-\frac{1}{%
x^{\alpha }}\right] ,
\end{equation*}%
so that $T=(\alpha a_{0}x_{0}^{\alpha })^{-1}$. By (\ref{Eq:PRA1}), 
\begin{equation}
P\left( T_{e}\leq t\right) \leq 1-\Phi \left( \frac{(\alpha
a_{0}t)^{-1/\alpha }-x_{0}}{\sqrt{T}}\right) ,  \label{Eq:AFE}
\end{equation}%
for every $0\leq t<T$. In order to estimate the value of $t$ for which $%
T_{e}\leq t$ with a probability of at most $0.05$, Equation (\ref{Eq:AFE})
yields 
\begin{equation*}
1-\Phi \left( \frac{(\alpha a_{0}t)^{-1/\alpha }-x_{0}}{\sqrt{T}}\right)
\leq 0.05
\end{equation*}%
whence it follows that $t=\frac{1}{\alpha a_{0}}[x_{0}+\sqrt{T}\Phi
^{-1}(0.95)]^{-\alpha }$. Proposition \ref{Prop:2} and monotonicity on the
integrand imply that 
\begin{equation*}
\begin{array}{lll}
\displaystyle{P\left( T_{e}\leq t\right) } & \leq & \displaystyle{P\left(
B(\infty )\leq \widetilde{A}(t)\right) } \\ 
& {\leq } & {P\left( \frac{1}{\alpha x_{0}^{\alpha }}\leq
\int_{0}^{t}a_{0}\left( \frac{1}{x_{0}}|\widehat{W}_{t}|+1\right) ^{1+\alpha
}ds\right) } \\ 
& \leq & \displaystyle{P\left( \frac{1}{\alpha x_{0}^{\alpha }}\leq
a_{0}t\left( \frac{1}{x_{0}}|\widehat{W}_{t}|+1\right) ^{1+\alpha }\right) }
\\ 
& {=} & {1-\Phi \left( \frac{x_{0}}{\sqrt{t}}\left( (\alpha
a_{0}x_{0}^{\alpha }t)^{-1/(1+\alpha )}-1\right) \right) }.%
\end{array}%
\end{equation*}%
This last estimate of $P(T_{e}\leq t)$ is better than that given by (\ref%
{Eq:AFE}), in view of the fact that $\alpha a_{0}x_{0}^{\alpha }t>1$.}
\end{example}

\subsubsection*{Acknowledgments}

M. J. Ceballos-Lira wishes to acknowledge the financial support of the
Mexican Council for Science and Technology (CONACYT) to pursue postgraduate
studies in Universidad Ju\'{a}rez Aut\'{o}noma de Tabasco (UJAT); he also
wishes to thank UJAT and the Universidad Aut\'{o}noma de Aguascalientes
(UAA) for additional, partial, financial support. J. Villa acknowledges the
partial support of CONACYT grant 118294, and grant PIM08-2 at UAA.

\end{document}